\documentclass[amstex,12pt, amssymb]{article}

\usepackage{mathtext}
\usepackage[cp1251]{inputenc}
\usepackage[T2A]{fontenc}
\usepackage[dvips]{graphicx}
\usepackage{amsmath}
\usepackage{amssymb}
\usepackage{amsxtra}
\usepackage{latexsym}
\usepackage{ifthen}

\newcounter{lemma}[section]

\newcounter{corollary}[section]

\newcounter{remark}[section]

\newcounter{theorem}[section]

\newcounter{proposition}[section]

\newcounter{example}

\numberwithin{equation}{section}

\begin{document}

\markboth{O.~DOVHOPIATYI, E.~SEVOST'YANOV}{\centerline{On boundary
H\"{o}lder logarithmic continuity ...}}

\def\cc{\setcounter{equation}{0}
\setcounter{figure}{0}\setcounter{table}{0}}

\overfullrule=0pt


\author{OLEKSANDR DOVHOPIATYI, EVGENY SEVOST'YANOV}

\title{
{\bf On boundary H\"{o}lder logarithmic continuity of mappings in
some do\-mains}}

\date{\today}
\maketitle

\begin{abstract}
We study mappings satisfying some estimate of distortion of modulus
of families of paths. Under some conditions on definition and mapped
domains, we have proved that these mappings are logarithmic
H\"{o}lder continuous at boundary points.
\end{abstract}


\section{Introduction}

This manuscript is devoted to the so-called inverse Poletsky
inequality established in many classes of mappings (see, e.g.,
\cite[Theorem~3.2]{MRV$_1$}, \cite[Theorem~6.7.II]{Ri} and
\cite[Theorem~8.5]{MRSY}). Recall that, mappings with a bounded
distortion as well as quasiconformal mappings satisfy the
inequalities
\begin{equation}\label{eq20}
M(\Gamma)\leqslant N(f, D)K_O(f)M(f(\Gamma))
\end{equation}
for any family of paths $\Gamma$ in $D,$ where $M$ denotes the
modulus of families of paths, $1\leqslant K_O(f)<\infty$ and $N(f,
D)$ is a maximal multiplicity of $f$ in $D$ (see~\cite{MRV$_1$}).
For more general classes of mappings, such an inequality has some
more general form (see below). Note that, the inequalities
$M(f(\Gamma))\leqslant K\cdot M(\Gamma),$ $1\leqslant K<\infty,$ are
very similar to~(\ref{eq20}) and were established by
Poletsky~\cite[Theorem~1, $\S\,4$]{Pol} for quasiregular mappings.
Precisely for this reason, we call the mappings in~(\ref{eq20})
mappings satisfying the inverse Poletsky inequality. However, we
will call this phrase mappings from a wider class.

\medskip
A Borel function $\rho:{\Bbb R}^n\,\rightarrow [0,\infty] $ is
called {\it an admissible} for a family $\Gamma$ of paths $\gamma$
in ${\Bbb R}^n,$ if the relation
\begin{equation}\label{eq1.4}
\int\limits_{\gamma}\rho (x)\, |dx|\geqslant 1
\end{equation}
holds for any locally rectifiable path $\gamma\in\Gamma.$
{\it A modulus} of $\Gamma $ is defined as follows:
\begin{equation}\label{eq1.3gl0}
M(\Gamma)=\inf\limits_{\rho \in \,{\rm adm}\,\Gamma}
\int\limits_{{\Bbb R}^n} \rho^n (x)\,dm(x)\,.
\end{equation}
Let $Q:{\Bbb R}^n\rightarrow [0, \infty]$ be a Lebesgue measurable
function. We say that {\it $f$ satisfies the inverse Poletsky
inequality}, if the relation
\begin{equation}\label{eq2*A}
M(\Gamma)\leqslant \int\limits_{f(D)} Q(y)\cdot\rho_*^n(y)\, dm(y)
\end{equation}
holds for any family of paths $\Gamma$ in $D$ and any $\rho_*\in
{\rm adm}\,f(\Gamma).$ Note that estimates~(\ref{eq2*A}) hold in
many classes of mappings (see, e.g., \cite[Theorem~3.2]{MRV$_1$},
\cite[Theorem~6.7.II]{Ri} and \cite[Theorem~8.5]{MRSY}).

\medskip
Recently we have proved logarithmic H\"{o}lder continuity for
mappings in~(\ref{eq2*A}) at the boundary of the unit ball
(see~\cite{Sev$_1$}). In this paper, we study similar mappings
between another type of domains. In particular, we deal with
mappings between the so-called quasiextremal distance domains and
convex domains. Note that, quasiextremal distance domains
(abbreviated $QED$-domains) are introduced by Gehring and Martio
in~\cite{GM} and are structures in which the modulus of the families
of paths is metrically related to the diameter of the sets. As for
the maps acting between these domains, here we obtained logarithmic
H\"{o}lder continuity at the corresponding points of the boundary.
In the next parts of the paper, we also study mappings of a domain
with a locally quasiconformal boundary (collared domains) onto a
convex domain. Besides that, we have considered mappings of some
regular domains which are defined as quasiconformal images of
domains with a locally quasiconformal boundary. For these classes of
mappings, we also have proved the H\"{o}lder logarithmic continuity
at the corresponding boundary points. It should be noted that for
the case of regular domains, this property is formulated in terms of
prime ends, and not in the Euclidean sense.

\medskip
A mapping $f:D\rightarrow {\Bbb R}^n$ is called {\it a discrete} if
$\{f^{-1}\left(y\right)\}$ consists of isolated points for any
$y\,\in\,{\Bbb R}^n,$ and {\it an open,} if the image of any open
set $U\subset D$ is an open set in ${\Bbb R}^n.$ A mapping $f$
between domains $D$ and $D^{\,\prime}$ is said to be {\it a closed}
if $f(E)$ is closed in $D^{\,\prime}$ for any closed set $E\subset
D$ (see, e.g., \cite[Section~3]{Vu$_1$}).

\medskip
In accordance with~\cite{GM}, a domain $D$ in ${\Bbb R}^n$ is called
a {\it domain of quasiextremal length}\index{domain of quasiextremal
length}\index{$QED$-domain} (a $QED$-{\it domain for short)} if
\begin{equation}\label{eq4***}
M(\Gamma(E, F, {\Bbb R}^n))\leqslant  A\cdot M(\Gamma(E, F, D))
\end{equation}
for some finite number $A\geqslant 1$ and all continua $E$ and $F$
in $D$. Observe that, a half-space or a ball are quasiextremal
distance domains, see~\cite[Lemma~4.3]{Vu$_2$}.

\medskip
Later, in the extended Euclidean space $\overline{{{\Bbb
R}}^n}={{\Bbb R}}^n\cup\{\infty\},$  we use the {\it spherical
(chordal) metric} $h(x,y)=|\pi(x)-\pi(y)|,$ where $\pi$ is a
stereographic projection of $\overline{{{\Bbb R}}^n}$ onto the
sphere $S^n(\frac{1}{2}e_{n+1},\frac{1}{2})$ in ${{\Bbb R}}^{n+1},$
and
$$h(x,\infty)=\frac{1}{\sqrt{1+{|x|}^2}}\,,$$
\begin{equation}\label{eq3C}
\ \ h(x,y)=\frac{|x-y|}{\sqrt{1+{|x|}^2} \sqrt{1+{|y|}^2}}\,, \ \
x\ne \infty\ne y
\end{equation}
(see, e.g., \cite[Definition~12.1]{Va}).
In what follows, given $A, B\subset \overline{{\Bbb R}^n}$ we set
$$h(A, B)=\inf\limits_{x\in A, y\in B}h(x, y)\,,\quad h(A)=\sup\limits_{x, y\in A}h(x ,y)\,,$$
where $h$ is a chordal metric in~(\ref{eq3C}). Similarly, we define
the Euclidean distance between sets and the Euclidean diameter by
the formulae
$$d(A, B)=\inf\limits_{x\in A, y\in B}|x-y|\,,\quad d(A)=\sup\limits_{x, y\in A}|x-y|\,.$$
Sometimes we also write ${\rm dist\,}(A, B)$ instead $d(A, B)$ and
${\rm\,diam\,}E$ instead $d(E),$ as well. As usually, we set
$$
B(x_0, r)=\{x\in {\Bbb R}^n: |x-x_0|<r\}\,,\qquad {\Bbb B}^n=B(0,
1)\,,$$
$$S(x_0,r) = \{ x\,\in\,{\Bbb R}^n : |x-x_0|=r\}\,.$$

\medskip
Given $\delta>0,$ domains $D, D^{\,\prime}\subset {\Bbb R}^n,$
$n\geqslant 2,$ a non-degenerate continuum $A\subset D^{\,\prime}$
and a Lebesgue measurable function $Q:D^{\,\prime}\rightarrow [0,
\infty]$ we denote by ${\frak S}_{\delta, A, Q }(D, D^{\,\prime})$ a
family of all open discrete and closed mappings $f$ of $D$ onto
$D^{\,\prime}$ satisfying the relation~(\ref{eq2*A}) such
that~$h(f^{\,-1}(A),
\partial D)\geqslant~\delta.$
The following statement holds.

\medskip
\begin{theorem}\label{th1}
{\sl\,Let $Q\in L^1(D^{\,\prime}),$ let $D$ be a quasiextremal
distance domain, and let $D^{\,\prime}$ be a convex bounded domain.
Then any $f\in {\frak S}_{\delta, A, Q}(D, D^{\,\prime})$ has a
continuous extension
$f:\overline{D}\rightarrow\overline{D^{\,\prime}}$ and, besides
that, for any $x_0\in\partial D$ there is its neighborhood $U$ and
constants $C_n=C(n, A, D, D^{\,\prime})>0$ and $r_0=r_0(x_0)>0$ such
that
\begin{equation}\label{eq2C}
|\overline{f}(x)-\overline{f}(y)|\leqslant\frac{C_n\cdot (\Vert
Q\Vert_1)^{1/n}}{\log^{1/n}\left(1+\frac{r_0}{|x-y|}\right)}
\end{equation}
for any $x, y\in U\cap \overline{D},$
where $\Vert Q\Vert_1$ is a norm of the function $Q$ in $L^1(D).$
 }
\end{theorem}

Consider the following definition that has been proposed
N\"akki~\cite{Na}, cf.~\cite{KR$_1$}. The boundary of a domain $D$
is called {\it locally quasiconformal,} if every point
$x_0\in\partial D$ has a neighborhood $U,$ for which there exists a
quasiconformal mapping $\varphi$ of $U$ onto the unit ball ${\Bbb
B}^n\subset{\Bbb R}^n$ such that $\varphi(\partial D\cap U)$ is the
intersection of the unit sphere ${\Bbb B}^ n$ with a coordinate
hyperplane $x_n=0,$ where $x=(x_1,\ldots, x_n.)$ Note that, with
slight differences in the definition, domains with such boundaries
are also called collared domains.

\medskip
\begin{theorem}\label{th3}
{\sl\,Let $Q\in L^1(D^{\,\prime}),$ let $D$ be a domain with a
locally quasiconformal boundary, and let $D^{\,\prime}$ be a bounded
convex domain. Then any $f\in {\frak S}_{\delta, A, Q}(D,
D^{\,\prime})$ has a continuous extension
$f:\overline{D}\rightarrow\overline{D^{\,\prime}},$ while, for any
$x_0\in\partial D$ there is a neighborhood $V$ of $x_0$ and some
numbers $C_n=C(n, A, D, D^{\,\prime}, x_0)>0,$ $r_0=r_0(n, A, x_0,
D)>0$ and $0<\alpha=\alpha(x_0)\geqslant 1$ such that
\begin{equation}\label{eq2C_1}
|\overline{f}(x)-\overline{f}(y)|\leqslant\frac{C_n\cdot (\Vert
Q\Vert_1)^{1/n}}{\log^{1/n}\left(1+\frac{r_0}{|x-y|^{\alpha}}\right)}
\end{equation}
for any $x, y\in V\cap \overline{D},$
where $\Vert Q\Vert_1$ is a norm of $Q$ in $L^1(D).$
 }
\end{theorem}

\medskip Using the result of the previous theorem, it is also
possible to obtain a statement about H\"{o}lder's logarithmic
continuity for bad boundaries in terms of prime ends.

The definition of a prime end used below may be found in~\cite{ISS},
cf.~\cite{KR$_1$}. We say that a bounded domain $D$ in ${\Bbb R}^n$
is {\it regular}, if $D$ can be quasiconformally mapped to a domain
with a locally quasiconformal boundary whose closure is a compact in
${\Bbb R}^n,$ and, besides that, every prime end in $D$ is regular.
Note that space $\overline{D}_P=D\cup E_D$ is metric, which can be
demonstrated as follows. If $g:D_0\rightarrow D$ is a quasiconformal
mapping of a domain $D_0$ with a locally quasiconformal boundary
onto some domain $D,$ then for $x, y\in \overline{D}_P$ we put:
\begin{equation}\label{eq1A}
\rho(x, y):=|g^{\,-1}(x)-g^{\,-1}(y)|\,,
\end{equation}
where the element $g^{\,-1}(x),$ $x\in E_D,$ is to be understood as
some (single) boundary point of the domain $D_0.$ The specified
boundary point is unique and well-defined, see
e.g.~\cite[Theorem~2.1, Remark~2.1]{IS}, cf.~\cite[Theorem~4.1]{Na}.

\medskip
The following statement holds.

\begin{theorem}\label{th4}
{\sl\,Let $Q\in L^1(D^{\,\prime}),$ let $D$ be regular domain, and
let $D^{\,\prime}$ be a bounded convex domain. Then any $f\in {\frak
S}_{\delta, A, Q}(D, D^{\,\prime})$ has a continuous extension
$f:\overline{D}_P\rightarrow\overline{D^{\,\prime}},$ in addition,
for any $P_0\in E_D$ there is a neighborhood $V$ of this point in
$(\overline{D}_P, \rho)$ and numbers $C_n=C(n, A, D, D^{\,\prime},
x_0)>0,$ $r_0=r_0(n, A, x_0, D)>0$ and
$0<\alpha=\alpha(x_0)\geqslant 1$ such that
\begin{equation}\label{eq2C_2}
|\overline{f}(P_1)-\overline{f}(P_2)|\leqslant\frac{C_n\cdot (\Vert
Q\Vert_1)^{1/n}}{\log^{1/n}\left(1+\frac{r_0}{\rho^{\alpha}(P_1,
P_2)}\right)}
\end{equation}
for any $P_1, P_2\in V,$
where $\Vert Q\Vert_1$ is a norm of $Q$ in $L^1(D).$
 }
\end{theorem}

\section{Auxiliary lemmas}
Before proving the main statements, we prove the following important
lemma, which is proved in~\cite[Lemma~2.1]{Sev$_1$} for the case of
the unit ball.

\medskip
\begin{lemma}\label{lem1}
{\sl\, Let $D$ and $D^{\,\prime}$ be domains satisfying the
conditions of Theorem~\ref{th1}, and let $E$ be a continuum
in~$D^{\,\prime},$ $Q\in L^1(D^{\,\prime}).$ Then there is
$\delta_1>0$ such that ${\frak S}_{\delta, A, Q }\subset {\frak
S}_{\delta_1, E, Q}.$ In other words, if $f$ is an open discrete and
closed mapping of $D$ onto $D^{\,\prime}$ satisfying the
condition~(\ref{eq2*A}) such that $h(f^{\,-1}(A),
\partial D)\geqslant~\delta,$ then there is $\delta_1>0,$ which does not
depend on $f,$ such that  $h(f^{\,-1}(E),
\partial D)\geqslant~\delta_1.$ }
\end{lemma}

\medskip
\begin{proof}
We will generally use the scheme of the
proof~\cite[Lemma~2.1]{Sev$_1$}. Let us prove Lemma~\ref{lem1} from
the opposite. Suppose that its conclusion is not true. Then, there
are sequences $y_m\in E,$ $f_m\subset{\frak S}_{\delta, A, Q }$ and
$x_m\in D$ such that $f_m(x_m)=y_m$ and $h (x_m,
\partial D)\rightarrow 0$ as $m\rightarrow\infty.$ Without loss of
generality, we may assume that $x_m\rightarrow x_0$ as
$m\rightarrow\infty,$ where $x_0$ may be equal to $\infty$ if $D$ is
unbounded. By Theorem~3.1 in~\cite{SSD}, it follows that $f_m$ has a
continuous extension to~$x_0,$ moreover, the family
$\{f_m\}_{m=1}^{\infty }$ is equicontinuous at $x_0$ (see, e.g.,
\cite[Theorem~1.2]{SSD}). Then, for any $\varepsilon>0$ there is
$m_0\in {\Bbb N}$ such that $h(f_m(x_m), f_m(x_0))<\varepsilon$ for
$m\geqslant m_0.$ On the other hand, since $f_m$ is closed,
$f_m(x_0)\in \partial D^{\,\prime}.$ Due to the compactness of the
space $\overline{{\Bbb R}^n}$ and the closure of $\partial
D^{\,\prime},$ we may assume that $f_m(x_0)$ converges to some $B\in
\partial D^{\,\prime}$ as $m\rightarrow\infty.$ Therefore, by the triangle inequality,
$$h(f_m(x_m), f_m(x_0))\geqslant h(f_m(x_m), B)-h(B, f_m(x_0))\geqslant \frac{1}{2}
\cdot h(E, \partial D^{\,\prime})$$
for sufficiently large $m\in{\Bbb N}.$ Finally, we have a
contradiction: $h(f_m(x_m), f_m(x_0))\geqslant \delta_0,$
$\delta_0:=\frac{1}{2 } \cdot h(E, \partial D^{\,\prime})$ and, at
the same time, $h(f_m(x_m), f_m(x_0))<\varepsilon$ for $m\geqslant
m_0.$ The resulting contradiction refutes the original assumption.
The lemma is proved.~$\Box$
\end{proof}

\medskip
The following lemma was proved in the case where the domain
$D^{\,\prime}$ is the unit ball (see the proof of Theorem~1.1
in~\cite{Sev$_1$}). For an arbitrary convex domain, its proof is
significantly more difficult, since the previous methodology relied
significantly on the geometry of the ball.

\medskip
\begin{lemma}\label{lem2}
{\sl\, Let $D^{\,\prime}$ be a bounded convex domain in ${\Bbb
R}^n,$ $n\geqslant 2,$ and let $B(y_*, \delta_*/2)$ be a ball
centered at the point $y_*\in D^{\,\prime},$ where $\delta_*:=d(y_*,
\partial D^{\,\prime}).$ Let $z_0\in \partial D^{\,\prime}.$ Then
for any points $A , B\in B(z_0, \delta_*/8)\cap D^{\,\prime}$ there
are points $C, D\in \overline{B(y_*, \delta_*/2)},$ for which the
segments $[A, C]$ and $[B, D]$ are such that
\begin{equation}\label{eq13}
{\rm dist\,}([A, C], [B, D])\geqslant C_0\cdot |A-C|\,,
\end{equation}
where $C_0>0$ is some constant that depends only on $\delta_*$ and
$d(D^{\,\prime}).$
  }
\end{lemma}

\medskip
\begin{proof}
We put
$$\varepsilon_0:=|A-B|<\delta_0:=\delta_*/4\,.$$
Let us join the points $A$ and $y_*$ with the segment $I.$ Due to
the convexity of $D^{\,\prime},$ the segment $I$ completely belongs
to $D^{\,\prime}.$ Points $A,$ $y_*$ and $B$ form the plane $P$ (if
$A$, $y_*$ and $B$ belong to some straight line, then we define by
$P$ any plane, consisting $A$ and $B$). Consider the circle on this
plane
$$S=\{z\in P:
|z-A|=\varepsilon_0\}\,.$$
The position of the point $z=B$ on the circle $S$ is completely
determined by the angle $\varphi,$ $\varphi\in [-\pi, \pi)$ between
the vector $y_*-A$ and the radius vector $B-A$ of the point $z.$
Points on the circle are further denoted by polar coordinates using
pairs $z=(\varepsilon_0, \varphi).$ Our the further goal is to
investigate the main three cases regarding the intervals of change
of this angle.

\medskip
{\bf Case 1.} ''Large angles'': $\varphi\in [\pi/4, 3\pi/4],$ or
$\varphi\in [-\pi/4, -3\pi/4],$ see Figure~\ref{fig1A} for
illustrations.
\begin{figure}[h]
\centerline{\includegraphics[scale=0.5]{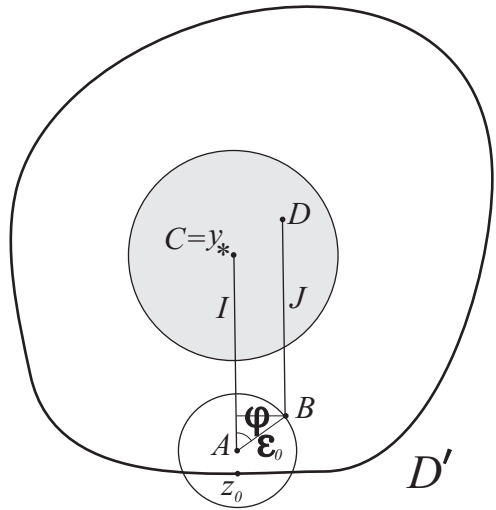}} \caption{To
proof of {\bf Theorem 1}, case~{\textbf{1}}}\label{fig1A}
\end{figure}
Put $C:=y_*.$ Let us limit ourselves to considering the case
$\varphi\in [\pi/4, 3\pi/4]$ (a case $\varphi\in [-\pi/4, -3\pi/4]$
may be considered similarly). Consider the ray
$$r=r(t)=B+te\,,\quad
e=(y_*-A)/|y_*-A|\,,\quad t>0.\,$$
By construction, the ray $r$ is parallel to the segment $I.$ When
$t=|y_*-A|,$ we have $r(|y_*-A|)=B+y_*-A$ and
$|r(|y_*-A|)-y_*|=|B-A|=\varepsilon_0<\delta_0,$ that is, point
$r(|y_*-A|)$ belongs to $E.$ Let $J$ be a segment of the ray $r,$
which is contained between the points $B$ and $D:=r(|y_*-A|)$ (it
belongs entirely to $D^{\,\prime}$ due to its convexity). The
distance between $I$ and $J$ is calculated as follows:
\begin{equation}\label{eq1}
{\rm dist}\,(I, J)=\varepsilon_0\sin\varphi\geqslant
\frac{\sqrt{2}}{2}\varepsilon_0\,.
\end{equation}
The consideration of the first case is completed, since we take $I$
and $J$ as segments $[A, C]$ and $[B, D],$ in addition, we put
$C_0:=\frac{\sqrt{2}}{2}$ in ~(\ref{eq13}).

\medskip
{\bf Case~2.} ''Small angles'': $\varphi\in [-\pi/4, \pi/4).$
Consider the case $\varphi\in [0, \pi/4)$ (the case $\varphi\in
[-\pi/4, 0)$ may be considered similarly). Let us draw a line
through the point $y_*$, which belongs to the plane $P$ and is
orthogonal to the vector $y_*-A,$ see Figure~\ref{fig2A}.
\begin{figure}[h]
\centerline{\includegraphics[scale=0.5]{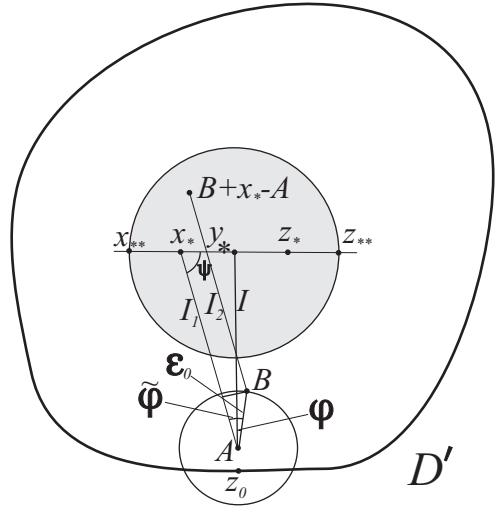}} \caption{To
proof of {\bf Theorem 1}, case~{\textbf{2}}}\label{fig2A}
\end{figure}
This line has exactly two points of intersection with the boundary
of the circle~$\overline{B(y_*, \delta_*/2)}\cap P.$ We denote these
points by $x_{**}$ and $z_{**}.$ Let $x_{**}$ be a point such that
the angle between by vectors $y_*-A$ and $x_{**}-A$ is negative. Let
us put
$$x_*=\frac{x_{**}+y_*}{2}\,,\qquad z_*=\frac{z_{**}+y_*}{2}\,.$$
Furthermore, let $\psi$ denote the angle $\angle Ax_*y_*.$ Then
\begin{equation}\label{eq6}
\tan
\psi=\frac{|A-y_*|}{\frac{\delta_*}{4}}<\frac{4d(D^{\,\prime})}{\delta_*}\,.
\end{equation}
In addition,
$$\tan
\psi=\frac{|A-y_*|}{\frac{\delta_*}{4}}\geqslant
\frac{\frac{\delta_*}{2}}{\frac{\delta_*}{4}}= 2\,,$$
from here
\begin{equation}\label{eq4A}
\psi\geqslant \arctan 2>\frac{\pi}{4}\,, \qquad
\frac{\pi}{2}-\psi\leqslant \frac{\pi}{2}-\arctan 2<\frac{\pi}{4}\,.
\end{equation}
Accordingly, the angle $\angle y_*Ax_*$ is equal to
$\frac{\pi}{2}-\psi.$ Let
$$\widetilde{\varphi}:=\varphi+\frac{\pi}{2}-\psi\,,$$
where $\psi$ is defined as above. Due to~(\ref{eq4A}), since
$\varphi\in [0, \pi/4)\,,$
\begin{equation}\label{eq5}
0<\widetilde{\varphi}\leqslant\frac{\pi}{4}+\frac{\pi}{2}-\arctan
2=\frac{3\pi}{4}-\arctan 2<\frac{\pi}{2}\,.
\end{equation}
Consider the ray $r(t)=B+t\frac{x_*-A}{|x_*-A|},$ $t>0.$ If
$t=|x_*-A|,$ then
$$r(|x_*-A|)=B+x_*-A\,.$$
Note that, $r(|x_*-A|)\in \overline{B\left(y_*,
\frac{\delta_*}{2}\right)}.$ Indeed, by the triangle inequality,
$$
|r(|x_*-A|)-y_*|=$$
\begin{equation}\label{eq2A}
=\bigl|B+x_*-A-y_*\bigr|\leqslant |x_*-y_*|+|B-A|\leqslant
\frac{\delta_*}{4}+\frac{\delta_*}{4}= \frac{\delta_*}{2}\,.
\end{equation}
Now let $I_1$ be the segment joining the points $A$ and $C:=x_{*},$
and let $I_2$ be the segment joining the points $B$ and
$D:=r(|x_*-A|)=B+x_*-A.$ By the construction and due
to~(\ref{eq2A}),
\begin{equation}\label{eq3A}
I_1\cap \overline{B\left(y_*,
\frac{\delta_*}{2}\right)}\ne\varnothing\ne I_2\cap
\overline{B\left(y_*, \frac{\delta_*}{2}\right)}\,.
\end{equation}
Note that
\begin{equation}\label{eq2}
{\rm dist}\,(I_1, I_2)=\varepsilon_0\sin \widetilde{\varphi}\,.
\end{equation}
Given the formulae~(\ref{eq6}) and~(\ref{eq5}), we obtain that
$$\tan\widetilde{\varphi}=\tan\left(\varphi+\frac{\pi}{2}-\psi\right)\geqslant$$
\begin{equation}\label{eq7}
\geqslant \tan\left(\frac{\pi}{2}-\psi\right)=\cot\psi \geqslant
\frac{\delta_*}{4d(D^{\,\prime})}\,.
\end{equation}
Due to~(\ref{eq7}) and taking into account~(\ref{eq5}), we obtain
that
\begin{equation}\label{eq8A}\sin\widetilde{\varphi}=\tan\widetilde{\varphi}\cdot
\cos\widetilde{\varphi}\geqslant
\frac{\delta_*}{4d(D^{\,\prime})}\cdot\cos\left(\frac{3\pi}{4}-\arctan
2\right)>0\,.
\end{equation}

\medskip
Now, by~(\ref{eq2}) and~(\ref{eq8A}), we obtain that
\begin{equation}\label{eq9}
{\rm dist}\,(I_1,
I_2)\geqslant\varepsilon_0\cdot\frac{\delta_*}{4d(D^{\,\prime})}\cdot\cos\left(\frac{3\pi}{4}-\arctan
2\right)>0\,.
\end{equation}
As we have already noted, the case $\varphi\in [-\pi/4, 0)$ may be
considered similarly, namely, the construction of the proof is the
same ''up to symmetry'' with respect to the segment joining the
points $A $ and $y_*.$ In particular, the role of the point $x_*$
will be played by the point $z_*;$ the role of the segment joining
the points $A$ and $x_*$ preforms the segment joining the points $A$
and $z_*$, etc.

\medskip
The inequality~(\ref{eq9}) completes the consideration of
case~\textbf{2}, because we may put $[A, C]:=I_1$ and $[B, D]:=I_2,$
in addition, we set
$C_0:=\frac{\delta_*}{4d(D^{\,\prime})}\cdot\cos\left
(\frac{3\pi}{4}-\arctan 2\right)$ in the inequality~(\ref{eq13}).

\medskip
\medskip
{\bf Case~3.} ''Very large angles:'' or $\varphi\in (3\pi/4, \pi],$
or $\varphi\in (-\pi, -3\pi/4).$ This case is very similar in
case~\textbf{2} at its discretion. It is enough to consider the
situation $\varphi\in (3\pi/4, \pi],$ because the second situation
$\varphi\in (-\pi, -3\pi/4)$ is considered similarly and differs
from the first one only by a certain ''symmetry'' relative to the
line passing through the points $A$ and $y_*.$ Let the points $ z_*$
and $z_{**}$ are exactly as defined in case~\textbf{2}. Consider the
segment $J_1,$ joining the points $A$ and $C:=z_*$ (it completely
belongs to $D^{\,\prime}$ due to its convexity). Consider a ray
\begin{equation}\label{eq10}
l(t)=B+t\cdot\frac{z_*-A}{|z_*-A|},\qquad t>0\,.
\end{equation}
For $t=|z_*-A|$, we obtain the point $D:=l(|z_*-A|)=B+z_*-A,$ see
Figure~\ref{fig3A}.
\begin{figure}[h]
\centerline{\includegraphics[scale=0.5]{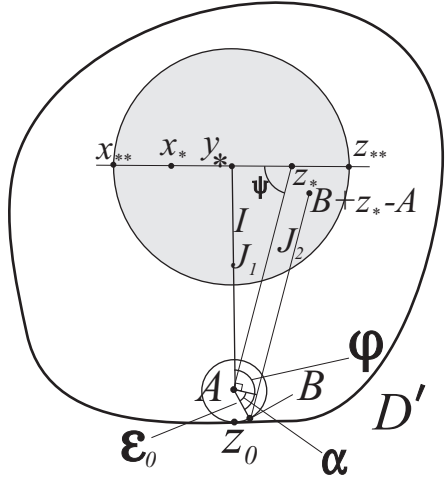}} \caption{To
proof of {\bf Theorem 1}, case~{\textbf{3}}}\label{fig3A}
\end{figure}
Let us prove that this point belongs to the ball
$\overline{B\left(y_*, \frac{\delta_*}{2}\right)}.$ Indeed, by the
triangle inequality
$$
|l(|z_*-A|)-y_*|=$$
\begin{equation}\label{eq2A_2}
=\bigl|B+z_*-A-y_*\bigr|\leqslant |z_*-y_*|+|B-A|\leqslant
\frac{\delta_*}{4}+\frac{\delta_*}{4}= \frac{\delta_*}{2}\,.
\end{equation}
Let $J_2$ be the segment joining the points $B$ and $D:=B+z_*-A.$ By
the construction and due to~(\ref{eq2A_2}),
\begin{equation}\label{eq3A_2}
J_1\cap \overline{B\left(y_*,
\frac{\delta_*}{2}\right)}\ne\varnothing\ne J_2\cap
\overline{B\left(y_*, \frac{\delta_*}{2}\right)}\,.
\end{equation}
Note that
\begin{equation}\label{eq2_2}
{\rm dist}\,(J_1, J_2)=\varepsilon_0\cos \alpha\,,
\end{equation}
where
$$\alpha=\varphi-\frac{\pi}{2}-\left(\frac{\pi}{2}-\psi\right)=
\varphi+\psi-\pi$$
and $\psi=\angle y_*z_*A.$ Now we need to find the lower bound for
$\cos \alpha.$

\medskip
From the triangle $\vartriangle y_*z_*A$ we obtain that
$$1\leqslant\tan \psi=\frac{|A-y_*|}{|y_*-z_*|}
\leqslant\frac{d(D^{\,\prime})}{\frac{\delta_*}{4}}=
\frac{4d(D^{\,\prime})}{\delta_*}\,,$$
from here $\frac{\pi}{4}\leqslant\psi\leqslant \arctan
\frac{4d(D^{\,\prime})}{\delta_*}<\frac{\pi}{2}.$ Then
\begin{equation}\label{eq11}
0=\frac{3\pi}{4}+\frac{\pi}{4}-\pi\leqslant \alpha:= \varphi+
\psi-\pi\leqslant \pi+\arctan
\frac{4d(D^{\,\prime})}{\delta_*}-\pi=\arctan
\frac{4d(D^{\,\prime})}{\delta_*}<\frac{\pi}{2}\,.
\end{equation}
By~(\ref{eq2_2}) and~(\ref{eq11}), we obtain that
\begin{equation}\label{eq12}
{\rm dist}\,(J_1, J_2)=\varepsilon_0\cos \arctan
\frac{4d(D^{\,\prime})}{\delta_*}\,.
\end{equation}
Actually, inequality~(\ref{eq12}) completes the consideration of
case~\textbf{3}, since we may put $[A, C]:=J_1$ and $[B, D]:=J_2,$
in addition, we put $C_0:=\cos \arctan
\frac{4d(D^{\,\prime})}{\delta_*}$ in~(\ref{eq13}).

\medskip
Finally, for all cases~\textbf{1}, \textbf{2} and \textbf{3}
in~(\ref{eq13}) we may put
$$C_0=\min\left\{\frac{\sqrt{2}}{2},
\frac{\delta_*}{4d(D^{\,\prime})}\cdot\cos\left(\frac{3\pi}{4}-\arctan
2\right), \cos \arctan
\frac{4d(D^{\,\prime})}{\delta_*}\right\}\,.$$
\end{proof}

\section{Proof of Theorem~\ref{th1}}

The possibility of a continuous extension of the mapping $f$ to the
boundary of the domain~$D$ follows by Theorem~3.1 in~\cite{SSD}. In
particular, the weak flatness of $\partial D$ is a consequence of
the fact that $D$ is a $QED$ domain (see, e.g.,
\cite[Lemma~2]{SevSkv}), in addition, any convex domain is locally
connected at its boundary because its intersection with the ball
centered at the boundary point is also a convex set.

\medskip
Let us prove logarithmic H\"{o}lder continuity~(\ref{eq2C}).
Put~$x_0 \in \partial D,$ and let $y_*\in D^{\,\prime}$ be an
arbitrary point of the domain $D^{\,\prime}.$ We may assume that
$x_0\ne \infty.$ Put $\delta_*:=d(y_*, \partial D^{\,\prime}).$ Let
$E=\overline{B(y_*, \delta_*/2)}\subset D^{\,\prime}.$ By
Lemma~\ref{lem1} there exists $\delta_1>0$ such that $h(f^{\,-1}(E),
\partial D)\geqslant \delta_1$ for all $f\in {\frak S}_{\delta, A,
Q}.$ In addition, due to Theorem~1.2 in~\cite{SSD} the
family~${\frak S}_{\delta, A, Q}$ is equicontinuous
in~$\overline{D},$ for a number $\delta_*/8$ there is a neighborhood
$U\subset B(x_0, \delta_1/2)$ of $x_0$ such that
$|f(x)-f(x_0)|<\delta_*/8$ for all $x, x_0\in U\cap D$ and all
$f\in{\frak S}_{\delta, A, Q}.$ Let $x, y\in U\cap D$ and
$$\varepsilon_0:=|f(x)-f(y)|<\delta_0:=\delta_*/4\,.$$
Lets us apply Lemma~\ref{lem2} for points $A=f(x),$ $B=f(y)$ and
$z_0=f(x_0).$ According to this lemma, there are segments $I\supset
A$ and $J\supset B$ in $D^{\,\prime}$ such that $I\cap
E\ne\varnothing\ne J\cap E,$ moreover
\begin{equation}\label{eq14}
{\rm dist\,}(I, J)\geqslant C_0\cdot |f(x)-f(y)|\,,
\end{equation}
where the constant $C_0$ depends only on the continuum $E$ and the
domain $D^{\,\prime}.$

\begin{figure}[h]
\centerline{\includegraphics[scale=0.5]{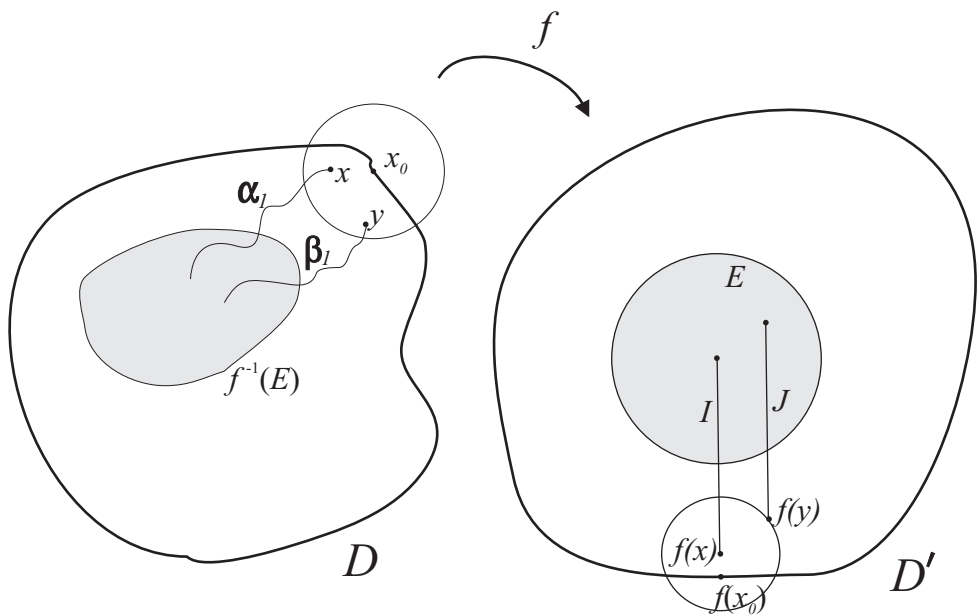}} \caption{To
proving the theorem~\ref{th1}}\label{fig1}
\end{figure}

\medskip
Let $\alpha_1,$ $\beta_1$ be total $f$-liftings of paths $I$ and $J$
starting at the points $x$ and $y,$, respectively (they exist due
to~\cite[Lemma~3.7]{Vu$_1$}, see Figure~\ref{fig1}).
By the definition, $|\alpha_1|\cap f^{\,-1}(E)\ne\varnothing\ne
|\beta_1|\cap f^{\,-1}(E).$ Since $h(f^{\,-1}(E),
\partial D)\geqslant \delta_1$ and $x, y\in B(x_0, \delta_1/2),$ then
\begin{equation}\label{eq4}
d(\alpha_1)\geqslant \delta_1/2\,,\quad d(\beta_1)\geqslant
\delta_1/2\,.
\end{equation}
Let
$$\Gamma:=\Gamma(\alpha_1, \beta_1, D)\,.$$
Then, on the one hand, by the inequality~(\ref{eq4***})
\begin{equation}\label{eq7A}
M(\Gamma)\geqslant (1/A_0)\cdot M(\Gamma(\alpha_1, \beta_1, {\Bbb
R}^n))\,,
\end{equation}
and on the other hand, by~\cite[Lemma~7.38]{Vu$_3$}
\begin{equation}\label{eq7B}
M(\Gamma(\alpha_1, \beta_1, {\Bbb R}^n))\geqslant
c_n\cdot\log\left(1+\frac1m\right)\,,
\end{equation}
where $c_n>0$ is some constant that depends only on $n,$
$$m=\frac{{\rm dist}(\alpha_1, \beta_1)}{\min\{{\rm diam\,}(\alpha_1),
{\rm diam\,}(\beta_1)\}}\,.$$
Then combining~(\ref{eq4}), (\ref{eq7A}) and~(\ref{eq7B}) and given
that ${\rm dist}\,(\alpha_1, \beta_1)\leqslant |x-y|,$ we obtain
that
\begin{equation}\label{eq7C}
M(\Gamma)\geqslant \widetilde{c_n}\cdot
\log\left(1+\frac{\delta_1}{2{\rm dist}(\alpha_1,
\beta_1)}\right)\geqslant \widetilde{c_n}\cdot
\log\left(1+\frac{\delta_1}{2|x-y|}\right)\,,
\end{equation}
where $\widetilde{c_n}>0$ is some constant depending only on $n$
and~$A_0$ from the definition of $QED$-domain.

\medskip
Now, let us to find some upper estimate for $M(\Gamma).$ Put
$$\rho(y)= \left\{
\begin{array}{rr}
\frac{1}{C_0\varepsilon_0}, & y\in D^{\,\prime},\\
0, & y\not\in D^{\,\prime}\,.
\end{array}
\right. $$
Note that, $\rho$ satisfies the relation~(\ref{eq1.4}) for
$f(\Gamma)$ by virtue of~(\ref{eq13}). By the definition of the
family ${\frak S}_{\delta, A, Q }$, we obtain that
\begin{equation}\label{eq14***}
M(\Gamma)\leqslant
\frac{1}{C^n_0\varepsilon_0^n}\int\limits_{D^{\,\prime}}
Q(y)\,dm(y)=C^{\,-n}_0\cdot \frac{\Vert
Q\Vert_1}{{|f(x)-f(y)|}^n}\,.
\end{equation}
By~(\ref{eq7C}) and (\ref{eq14***}), it follows that
$$\widetilde{c_n}\cdot \log\left(1+\frac{\delta_1}{2|x-y|}\right)\leqslant
C^{\,-n}_0\cdot\frac{\Vert Q\Vert_1}{{|f(x)-f(y)|}^n}\,.$$
The desired inequality ~(\ref{eq2C}) follows from the last relation,
where $C_n:=C^{\,-n}_0\cdot\widetilde{c_n}^{-1/n}$
and~$r_0=\delta_1/2.$

We have proved Theorem~\ref{th1} for the inner points $x, y\in U\cap
D.$ For the points $x, y\in U\cap \overline{D}$, this statement
follows by means passing to the limit $\overline{x}\rightarrow x$
and $\overline{y}\rightarrow y,$ $\overline{x}, \overline{y}\in
D.$~$\Box$

\medskip
The analogue of Theorem~\ref{th1} is also valid for mappings with a
fixed point in~$D.$ In order to formulate and prove the
corresponding statement, let us introduce the following definition.
For $a, b\in D$ and a Lebesgue measurable function
$Q:D^{\,\prime}\rightarrow [0, \infty]$ we denote by ${\frak F}_{a,
b, Q }$ the family of all open discrete and closed mappings $f$ of
the domain $D$ onto $D^{\,\prime}$ such that~$f(a)=b.$ The following
is true.

\medskip
\begin{theorem}\label{th2}
{\sl\,Let $Q\in L^1(D^{\,\prime}).$ Then any mapping $f\in {\frak
S}_{\delta, A, Q }$ has a continuous extension to a mapping
$f:\overline{D}\rightarrow\overline{D^{\,\prime}},$ while, for any
point $x_0\in\partial D$ there is a neighborhood $U$ of this point
and numbers $C_n=C(n, A, D, D^{\,\prime})>0$ and $r_0=r_0(x_0)>0$
such that the relation~(\ref{eq2C}) is fulfilled. }
\end{theorem}

\medskip
\begin{proof}
The possibility of a continuous extension of $f$ to $\partial D$
follows by~\cite[Theorem~3.1]{SSD}. Let us prove the logarithmic
H\"{o}lder continuity of the cooresponding family of extended
mappings.

\medskip
Put $E=\overline{B(b, \varepsilon_*)},$ where $\varepsilon_*={\rm
dist\,}(b, \partial D^{\,\prime}).$ Two cases are possible:

\medskip
1) there exists $\delta>0$ such that ${\rm dist}\,(f^{\,-1}(E),
\partial D)\geqslant\delta $ for all $f\in {\frak S}_{\delta, A, Q
}.$ In in this case, the desired statement follows by
Theorem~\ref{th1};

\medskip
2) There are sequences $f_m\in {\frak S}_{\delta, A, Q }$ and
$x_m\in D, y_m\in D^{\,\prime},$ $m=1,2,\ldots ,$ such that
$f_m(x_m)=y_m,$ $y_m\in E$ and ${\rm dist\,}(x_m,
\partial D)\rightarrow 0$ as $m\rightarrow\infty.$
Reasoning in the same way as in the proof of Lemma~\ref{lem1}, we
come to the conclusion that the family~${\frak S} _{\delta, A, Q }$
is not equicontinuous at at least one point $x_0\in \partial D,$
however, this contradicts the statement of Theorem~7.1
in~\cite{SSD}.

\medskip
Thus, the case 2) is impossible, and in the case~1) we have the
desired assertion of the theorem.~$\Box$
\end{proof}

\section{Proof of Theorem~\ref{th3}}

The possibility of a continuous extension of the mapping $f$ to the
boundary of $D$ follows by Theorem~3.1 in~\cite{SSD}. In particular,
locally quasiconformal boundaries of domains are weakly flat (see
\cite[Proposition~2.2]{ISS}, see also~\cite[Theorem~17.10]{Va}), and
convex domains are obviously locally connected at its boundary.

\medskip
Put~$x_0 \in \partial D.$ Let $y_*\in D^{\,\prime}$ be an arbitrary
point of $D^{\,\prime},$ $\delta_*:=d(y_*,
\partial D^{\,\prime})$ and $E=\overline{B(y_*, \delta_*/2)}\subset
D^{\,\prime}.$ By Lemma~\ref{lem1}, there exists $\delta_1>0$ such
that $h(f^{\,-1}(E),
\partial D)\geqslant \delta_1$ for all $f\in {\frak S}_{\delta, A,
Q}.$ Then $d(f^{\,-1}(E), \partial D)\geqslant \delta_1$ for any
$f\in {\frak S}_{\delta, A, Q} .$ In addition, since by Theorem~1.2
in~\cite{SSD} the family~${\frak S}_{\delta, A, Q}$ is
equicontinuous at~$\overline{D}, $ for $\delta_*/8$ there is a
neighborhood $U\subset B(x_0, \delta_1/4)$ of $x_0$ such that
$|f(x)-f(x_0)|<\delta_*/ 8$ for any $x, y\in U\cap D$ and all
$f\in{\frak S}_{\delta, A, Q}.$

\begin{figure}[h]
\centerline{\includegraphics[scale=0.5]{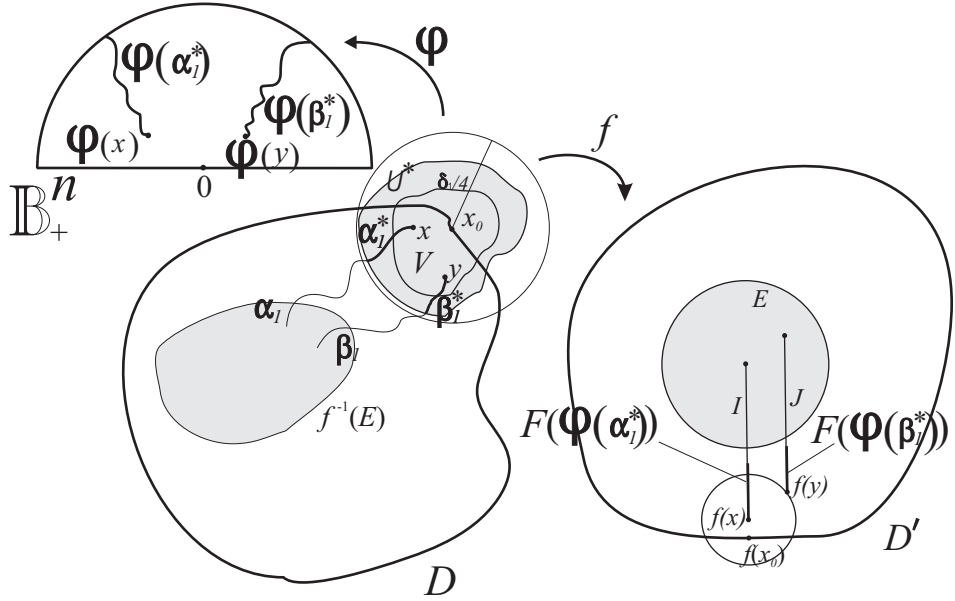}} \caption{To
the proof of Theorem~\ref{th3}}\label{fig4}
\end{figure}

By the definition of a locally quasiconformal boundary, there exists
a neighborhood $U^{\,*}$ of the point $x_0$ and a quasiconformal
mapping $\varphi:U^{\,*}\rightarrow {\Bbb B}^n,$
$\varphi(U^{\,*})={\Bbb B}^n,$ such that $\varphi(D\cap
U^{\,*})={\Bbb B}_{+}^n,$ where ${\Bbb B}_{+}^n=\{x\in {\Bbb B}^n:
x=(x_1,\ldots, x_n), x_n>0\}$ is a half-ball, see Figure~\ref{fig4}.
We may assume that $x_0\ne \infty$ and $\varphi(x_0)=0,$ and
$\overline{U^{\,*}}\subset U$ (see the proof of Theorem~17.10
in~\cite{Va}). Let $V$ be any neighborhood in $U^{\,*}$ such that
$\overline{V}\subset U^{\,*}.$ Let
\begin{equation}\label{eq23}
\delta_2:={\rm dist\,}(\partial V, \partial U^{\,*})\,.
\end{equation}
Consider the auxiliary mapping
\begin{equation}\label{eq16}
F(w):=f(\varphi^{\,-1}(w))\,,\qquad F:{\Bbb B}_{+}^n\rightarrow
U^{\,*}\,.
\end{equation}
Let $x, y\in V\cap D$ and
$$\varepsilon_0:=|f(x)-f(y)|<\delta_0:=\delta_*/4\,.$$
Apply Lemma~\ref{lem2} for points $A=f(x),$ $B=f(y)$ and
$z_0=f(x_0).$ According to this lemma, there are exist segments
$I\supset A$ and $J\supset B$ in $D^{\,\prime}$ such that $I\cap
E\ne\varnothing\ne J\cap E,$ moreover
\begin{equation}\label{eq14_1}
{\rm dist\,}(I, J)\geqslant C_0\cdot |f(x)-f(y)|\,,
\end{equation}
where the constant $C_0$ depends only on the continuum $E$ and the
domain $D^{\,\prime}.$

\medskip
Let $\alpha_1,$ $\beta_1$ be the whole $f$-lifts of the paths $I$
and $J$ with origins at the points $x$ and $y,$ respectively (they
exist by~\cite[Lemma~3.7]{Vu$_1$}). Then, by the definition,
$|\alpha_1|\cap f^{\,-1}(E)\ne\varnothing\ne |\beta_1|\cap
f^{\,-1}(E).$ Then
$$|\alpha_1|\cap U^{\*}\ne\varnothing\ne |\alpha_1|\cap({\Bbb R}^n\setminus
U^{\*})$$ and
$$|\beta_1|\cap U^{\*}\ne\varnothing\ne |\beta_1|\cap({\Bbb R}^n\setminus
U^{\*})\,.$$
Then, by~\cite[Theorem~1.I.5.46]{Ku}
\begin{equation}\label{eq21} |\alpha_1|\cap \partial U^{\*}\ne\varnothing\,,
|\beta_1|\cap
\partial U^{\*}\ne\varnothing\,.
\end{equation}
Similarly,
\begin{equation}\label{eq22} |\alpha_1|\cap \partial V\ne\varnothing\,,
|\beta_1|\cap
\partial V\ne\varnothing\,.
\end{equation}
Due to~(\ref{eq21}), $\alpha_1$ and $\beta_1$ contain subpaths
$\alpha^{\,*}_1$ and $\beta^{\,*}_1$ with origins at the points $x$
and $y$ which belong entirely in $U^{\,*}$ and have end points at
$\partial U^{\,*}.$ Due to~(\ref{eq23}), (\ref{eq21}) and
(\ref{eq22})
\begin{equation}\label{eq4_1}
d(\alpha^{\,*}_1)\geqslant \delta_2\,,\quad
d(\beta^{\,*}_1)\geqslant \delta_2\,.
\end{equation}
Consider the paths $\varphi(\alpha^{\,*}_1)$ and
$\varphi(\beta^{\,*}_1).$ Since $\varphi$ is quasiconformal mapping,
so $\varphi^{\,-1}$ is also quasiconformal. Thus, $\varphi^{\,-1}$
is H\"{o}lder continuous with some constant $\widetilde{C}>0$ and
some exponent $0<\alpha\leqslant 1$
(see~\cite[Theorem~1.11.III]{Ri}).
Let $\overline{x}, \overline{y}\in U^{\,*}$ be such that
$d(\alpha^{\,*}_1)=|\overline{x}-\overline{y}|.$ We put
$x^{\,*}=\varphi(\overline{x})$ and $y^{\,*}=\varphi(\overline{y}).$
Then
$$|x^{\,*}-y^{\,*}|^{\alpha}\geqslant\frac{1}{\widetilde{C}}\cdot
|\overline{x}-\overline{y}|=d(\alpha^{\,*}_1)\geqslant
\frac{1}{\widetilde{C}}\delta_2\,,$$
or
\begin{equation}\label{eq15}
|x^{\,*}-y^{\,*}|\geqslant
\left(\frac{1}{\widetilde{C}}\delta_2\right)^{1/\alpha}\,.
\end{equation}
From~(\ref{eq15}), we obtain that
$d(\varphi(\alpha^{\,*}_1))\geqslant
\left(\frac{1}{\widetilde{C}}\delta_2\right)^{1/\alpha}.$ Similarly,
$d(\varphi(\beta^{\,*}_1))\geqslant
\left(\frac{1}{\widetilde{C}}\delta_2\right)^{1/\alpha}.$
Let
$$\Gamma:=\Gamma(\varphi(\alpha^{\,*}_1), \varphi(\beta^{\,*}_1), {\Bbb B}_+^n)\,. $$
Note that, ${\Bbb B}_+^n$ is a bounded convex domain, so it is a
John domain (see~\cite[Remark~2.4]{MS}), hence it is a uniform
domain (see~\cite[Remark~2.13(c)]{MS}), therefore it is also a
$QED$-domain with some $A_0^{\,*}<\infty$ in~(\ref{eq4***})
(see~\cite[Lemma~2.18]{GM}). Then, on the one hand,
by~(\ref{eq4***})
\begin{equation}\label{eq7A_1}
M(\Gamma)\geqslant (1/A^{\,*}_0)\cdot
M(\Gamma(\varphi(\alpha^{\,*}_1), \varphi(\beta^{\,*}_1), {\Bbb
R}^n))\,,
\end{equation}
and on the other hand, by~\cite[Lemma~7.38]{Vu$_3$}
\begin{equation}\label{eq7B_1}
M(\Gamma(\varphi(\alpha^{\,*}_1), \varphi(\beta^{\,*}_1), {\Bbb
R}^n))\geqslant c_n\cdot\log\left(1+\frac1m\right)\,,
\end{equation}
where $c_n>0$ is some constant that depends only on $n,$
$$m=\frac{{\rm dist}(\varphi(\alpha^{\,*}_1), \varphi(\beta^{\,*}_1))}
{\min\{{\rm diam\,}(\varphi(\alpha^{\,*}_1)), {\rm
diam\,}(\varphi(\beta^{\,*}_1))\}}\,.$$
Then, combining~(\ref{eq7A_1}) and~(\ref{eq7B_1}) and taking into
account that ${\rm dist}\,(\varphi(\alpha^{\,*}_1),
\varphi(\beta^{\,*}_1))\leqslant |\varphi(x)-\varphi(y )|,$ we
obtain that
$$M(\Gamma)\geqslant \widetilde{c_n}\cdot \log\left
(1+\frac{\delta_2^{1/\alpha}}{(\widetilde{C})^{1/\alpha}{\rm
dist}(\alpha_1, \beta_1)}\right)\geqslant$$
\begin{equation}\label{eq7C_1}
\geqslant \widetilde{c_n}\cdot \log\left
(1+\frac{\delta_2^{1/\alpha}}{(\widetilde{C})^{1/\alpha}|\varphi(x)-\varphi(y)|}\right)\,,
\end{equation}
where $\widetilde{c_n}>0$ is some constant that depends only on $n$
and $A^{\,*}_0$ from the definition of $QED$-domain.

\medskip
Let us now establish an upper bound for $M(\Gamma).$ Note that, $F$
in~(\ref{eq16}) satisfies the relation~(\ref{eq2*A}) with the
function $\widetilde{Q}(x)=K_0\cdot Q(x)$ instead of $Q,$ where
$K_0\geqslant 1$ is the constant of a quasiconformality
of~$\varphi^{\,-1}.$ Let us put
$$\rho(y)= \left\{
\begin{array}{rr}
\frac{1}{C_0\varepsilon_0}, & y\in D^{\,\prime},\\
0, & y\not\in D^{\,\prime}\,,
\end{array}
\right. $$
where $C_0$ is the universal constant in inequality~(\ref{eq13}).
Note that $\rho$ satisfies the relation~(\ref{eq1.4}) for
$F(\Gamma)$ due to the relation~(\ref{eq13}). Then, by the
definition of ${\frak S}_{\delta, A, Q },$ due to the definition
of~$F$ in~(\ref{eq16}), we obtain that
\begin{equation}\label{eq14***_1}
M(\Gamma)\leqslant
\frac{1}{C^n_0\varepsilon_0^n}\int\limits_{D^{\,\prime}}
K_0Q(y)\,dm(y)=C^{\,-n}_0K_0\cdot \frac{\Vert
Q\Vert_1}{{|f(x)-f(y)|}^n}\,.
\end{equation}
It follows by~(\ref{eq7C_1}) and (\ref{eq14***_1}) that
$$\widetilde{c_n}\cdot \log\left
(1+\frac{\delta_2^{1/\alpha}}{(\widetilde{C})^{1/\alpha}|\varphi(x)-\varphi(y)|}\right)\
\leqslant C^{\,-n}_0K_0\cdot\frac{\Vert
Q\Vert_1}{{|f(x)-f(y)|}^n}\,.$$
From the last relation, by the H\"{o}lder continuity of the mapping
$\varphi$, it follows that
$$|f(x)-f(y)|\leqslant
C^{\,-n}_0{\widetilde{c_n}}^{-\frac{1}{n}}K_0\cdot\frac{\left(\Vert
Q\Vert_1\right)^{\frac{1}{n}}}{\log^{\frac{1}{n}}\left
(1+\frac{\delta^{1/\alpha}_2}{(\widetilde{C})^{1/\alpha}|\varphi(x)-\varphi(y)|}\right)
}\,\leqslant$$
$$\leqslant C^{\,-n}_0{\widetilde{c_n}}^{-\frac{1}{n}}K_0\cdot\frac{\left(\Vert
Q\Vert_1\right)^{\frac{1}{n}}}{\log^{\frac{1}{n}}\left
(1+\frac{\delta_2^{1/\alpha}}{
(\widetilde{C})^{(1/\alpha)+1}|x-y|^{\alpha }}\right) }\,,$$
which is the desired inequality~(\ref{eq2C_1}), where
$C_n:=C^{\,-n}_0\cdot\widetilde{c_n}^{-1/n}\cdot K_0$ and
$r_0=\frac{\delta_2^{1/\alpha}}{(\widetilde{C})^{1/\alpha+1}}.$

We proved Theorem~\ref{th3} for the inner points $x, y\in V\cap D.$
For $x, y\in V\cap \overline{D}$, this statement follows by means of
the transition to the limit $\overline{x}\rightarrow x$ and
$\overline{y}\rightarrow y,$ $\overline{x}, \overline{y}\in
D.$~$\Box$

\section{Proof of Theorem~\ref{th4}}

Let $f\in {\frak S}_{\delta, A, Q}(D, D^{\,\prime}).$ It is
sufficient to restrict ourselves to the case $P_1, P_2\in V\cap D.$
Since $D$ is a regular domain, there exists a quasiconformal mapping
$g^{\,-1}$ of a domain $D$ onto a domain $D_0$ with a locally
quasiconformal boundary, and, by the definition of the metric $\rho$
in~(\ref{eq1A}),
\begin{equation}\label{eq18}
\rho(P_1, P_2):=|g^{\,-1}(P_1)-g^{\,-1}(P_2)|\,. \
\end{equation}
Consider the auxiliary mapping
\begin{equation}\label{eq19}
F(x)=(f\circ g)(x)\,,\quad x\in D_0\,.
\end{equation}
Since $g^{\,-1}$ is quasiconformal, there is a constant $1\leqslant
K_1<\infty$ such that
\begin{equation}\label{eq17}
\frac{1}{K_1}\cdot M(\Gamma)\leqslant M(g^{\,-1}(\Gamma))\leqslant
K_1\cdot M(\Gamma)
\end{equation}
for any family of paths $\Gamma$ in $D_0.$ Considering
inequalities~(\ref{eq17}) and taking into account that $f$ satisfies
the relation~(\ref{eq2*A}), we obtain that also $F$ satisfies the
relation~(\ref{eq2*A}) with a new function
$\widetilde{Q}(x):=K_1\cdot Q(x).$ In addition, since $g$ is a fixed
homeomorphism, then $h(F^{\,- 1}(A), \partial
D)\geqslant~\delta_0>0,$ where $\delta_0>0$ is some fixed number.
Then Theorem~\ref{th2} may be applied to the map $F$. Applying this
theorem, we obtain that, for any point $x_0\in D_0$ there are $V$
neighborhood of this point and numbers $C^*_n=C(n, A, D_0, D^
{\,\prime})>0,$ $r_0=r_0(n, A, x_0, D_0)>0$ and
$0<\alpha=\alpha(x_0)\geqslant 1$ such that
\begin{equation}\label{eq2C_3}
|F(x)-F(y)|\leqslant\frac{C^*_nK_1^{\frac{1}{n}}\cdot (\Vert
Q\Vert_1)^{1/n}}{\log^{1/n}\left(1+\frac{r_0}{|x-y|^{\alpha}}\right)}
\end{equation}
for all $x, y\in V\cap D_0,$
where $\Vert Q\Vert_1$ is the norm of the function $Q$ in $L^1(D).$
Let $U:=g(V),$ $P_0:=g(x_0).$ Then, by definition, $U$ is a
neighborhood of the prime end $P_0\in E_D.$ If $P_1, P_2\in D_P\cap
U ,$ then $P_1=g(x)$ and $P_2=g(y)$ for some $x, y \in V\cap D_0.$
Taking into account the relation~(\ref{eq2C_3}) and using the
relation $|x-y|=|g^{\,-1}(P_1)-g^{\,-1}(P_2)|=\rho(P_1, P_2),$ we
obtain that
$$|F(g^{\,-1}(P_1))-F(g^{\,-1}(P_2))|\leqslant\frac{C^*_nK_1^{\frac{1}{ n}}\cdot (\Vert
Q\Vert_1)^{1/n}}{\log^{1/n}\left(1+\frac{r_0}{\rho^{\alpha}(P_1,
P_2)}\right)}\,,$$
or, due to~(\ref{eq19}),
$$|f(P_1)-f(P_2)|\leqslant\frac{C^*_nK_1^{\frac{1}{n}}\cdot (\Vert
Q\Vert_1)^{1/n}}{\log^{1/n}\left(1+\frac{r_0}{\rho^{\alpha}(P_1,
P_2)}\right)}\,.$$
The last ratio is desired if we put
here~$C_n:=C^*_nK_1^{\frac{1}{n}}.$~$\Box$

\medskip
{\bf \noindent Oleksandr Dovhopiatyi} \\
{\bf 1.} Zhytomyr Ivan Franko State University,  \\
40 Bol'shaya Berdichevskaya Str., 10 008  Zhytomyr, UKRAINE \\
alexdov1111111@gmail.com

\medskip
\medskip
{\bf \noindent Evgeny Sevost'yanov} \\
{\bf 1.} Zhytomyr Ivan Franko State University,  \\
40 Bol'shaya Berdichevskaya Str., 10 008  Zhytomyr, UKRAINE \\
{\bf 2.} Institute of Applied Mathematics and Mechanics\\
of NAS of Ukraine, \\
1 Dobrovol'skogo Str., 84 100 Slavyansk,  UKRAINE\\
esevostyanov2009@gmail.com

\end{document}